\documentclass[11pt]{amsart}

\usepackage[english]{babel}
\usepackage[latin1]{inputenc}
\usepackage{graphicx}
\usepackage{amsmath}
\usepackage{amssymb}
\usepackage{amscd}
\usepackage{color}

{\newtheorem{Theorem}{Theorem}[section]
\newtheorem{Lemma}[Theorem]{Lemma}
\newtheorem{Prop}[Theorem]{Proposition}

\newtheorem{Cor}[Theorem]{Corollary}

\newtheorem{Observ}[Theorem]{Remark}

\newtheorem{Conj}[Theorem]{Conjecture}

\setcounter{page}{1}

\def\Fol{{\f}(n,c)}

\def\j{\jmath}

\def\C{{\mathbb{C}}}
\def\K{{\mathbb{K}}}
\def\N{{\mathbb{N}}}
\def\P{{\mathbb{P}}}

\def\Z{{\mathbb{Z}}}

\def\f{{\mathcal{F}}}

\def\a{{\alpha}}
\def\si{{\sigma}}
\def\va{{\varphi}}

\def\o{\omega}

\def\Om{{\Omega}}

\def\tt{px\,dy-qy\,dx}
\def\rtt{x\,dy-y\,dx}

\def\O{{\mathcal{O}}}
\def\tcap{\pitchfork}
\title %[Foliations with radial ...]
{Foliations with a radial Kupka set on projective spaces}
\author[O. Calvo--Andrade]{Omegar Calvo--Andrade }
\thanks{MTM2004--07978, Spain.}
\dedicatory{}
\address{CIMAT: Ap. Postal 402, Guanajuato, 36000, Gto. M\'exico}
\email{omegar@cimat.mx}
\date{}
\begin{document}
\maketitle

\begin{abstract} We consider the set $K(n,c,\rtt)$ of
  codimension one holomorphic foliations on $\P^n,\,\, n\geq3$, with
  Chern class $c$, and with a compact, connected Kupka set of radial
  transversal type. We will prove that foliations in this set, have a
  rational first integral and define an irreducible component of the
  space of foliations.
\end{abstract}

\bigskip

\subsection{Introduction}
\label{sec:Int}

Let $n\geq 3$ and $c\geq 2$ be natural numbers. Consider a
differential 1--form in $\C^{n+1}$
\[ \o=a_0 dz_0+\cdots +a_n dz_n,\in   \]
where $a_j$ are homogeneous polynomials of degree $c-1$ in variables
$z_0,\dots,z_n$ with complex coefficients. Let us assume that
\[z_0a_0+\cdots + z_n a_n=0\]
so that $\o$ descends to the complex projective space $\P^n$ and
defines a global section of the twisted sheaf of 1--forms
$\Om^1_{\P^n}(c)$.

For a $\K$--vector space $V$, we denote $\P
V=V-\{\mathbf{0}\}/\K^{\ast}$ the projective space of linear
subspaces of $V$ and $\pi:V-\{\mathbf{0}\} \to \P V$ the quotient
map.

Consider the projective space $\P H^0(\P^n,\Om^1(c))$ and the
subset
\[\Fol :=\pi \left( \{\, \o\in
  H^0(\P^n,\Om^1_{\P^n}(c))-\mathbf{0}\, |\, \o\wedge d\o=0\, \}
 \right) \]
parameterizing 1--forms $\o$ such that they satisfy the Frobenius
integrability condition. This is the space of \emph{Chern class} $c$
\emph{foliations of codimension one on\/} $\P^n$. It is an algebraic
subset defined by quadratic equations, and has several irreducible
components. For instance, let us recall (see \cite{CL1},
 \cite{GML}) the following families of irreducible components:

The \emph{rational components} $R_n(a,b)\subset\Fol$ consisting of
integrable 1--forms of the type
  \[ \o=pg df-qfdg\in H^0(\P^n,\Om^1(c)) \]
where $c=a+b$ is a partition with $a,b$ natural numbers, $p,q$ are
the unique coprime numbers such that $pa=qb$ and $f,g$ are
homogeneous polynomials of respective degree $a,b$.

The leaves of the foliation $\o=pg df-qfdg\in R_n(a,b)$, are the
fibers of the rational map $\va=f^p/g^q$, and we say that $\va$ is
a \emph{rational first integral\/} of the foliation $\o$.

The \emph{singular set} $S(\o)=\{p\in\P^n | \,\o(p)=0  \}$, of a
foliation $\o\in \Fol,$ is the variety of zeros of the section $\o.$
We set $S_k(\o),$ the union of the irreducible components of $S(\o)$
of dimension $k$. It is always possible to suppose that
$codim(S(\o))\geq2$. Namely, if $S_{n-1}(\o)$ is not empty (i.e. if
$\o$ vanishes in codimension one). Then there exists a homogeneous
polynomial $F$ of maximal degree $0<e<c-1$ that divides $\o$. We
denote
\[ \overline{\o}=\frac{\o}{F}\in H^0(\P^n,\Om^1_{\P^n}(c-e)). \]
It is clear that $\overline{\o}$ is well defined up to a
multiplicative constant and it does not vanish in codimension one.

In this paper, the Rational Components $R_n(c/2,c/2)$ will be
characterized by its singular set, our main hypothesis, is that
$S_{n-2}(\o)$ is compact, connected and it is of Kupka type. The
Kupka set is defined by
\[ K(\o)= \{ p\in \P^n\,|\, \o(p)=0,\quad d\o(p)\neq0 \, \}\subset S(\o). \]

The main properties of Kupka sets, are the following
\begin{description}
\item[(a)] $K(\o)$ is smooth of codimension two.
\item[(b)] It is stable under deformations of the foliation.
\item[(c)] $K(\o)$ has \emph{local product structure\/}.
\end{description}

We are going to consider the subset of foliations
\[ K(n,c)   =  \{ \o\in\Fol | \, S_{n-2}=K(\o)
        \mbox{ is compact and connected } \} \]

The stability under deformations of the Kupka set, implies that $K(n,c)$
is an open subset of $\Fol$, its closure, is a union
of irreducible components. On the other hand, a fundamental fact about
a foliation $\o\in K(n,c)$, is that $\o$ has a rational first
integral and belongs to some irreducible rational component $R_n(a,b)$; if and
only if, its Kupka set $K(\o)$, is a complete
intersection \cite{CL}.
%. In this case $\o$ belongs to some rational component
% $R_n(a,b)$
Furthermore, it is conjectured that $K(\o)$ is always a complete
intersection \cite{CL}. There are already several partial positive
answers to this conjecture. In fact: in \cite[Theorem 3.5]{C} it
is proven that for foliations $\o\in K(n,c)$ which are not of
radial type, there always is a rational first integral, and in
\cite{B, B1}, we find the hypothesis $n\geq6$. The conjecture
remains open for the radial transversal type in small dimensions.

In this note we are going to prove the following result:

\begin{Theorem}\label{main} The set
$K(n,c,\rtt)\subset K(n,c)$ of foliations of radial type is not
empty if and only if $c$ is even, moreover there are homogeneous
polynomials $f,g$ of degree $c/2$ such that $K(\o)$ is the scheme
theoretically complete intersection  $K(\o)=\{f=g=0\}$ and
\[\o=f\,dg-g\,df. \]
\end{Theorem}

We observe that for any even natural number $c$, the generic element
of the irreducible component $R_n(c/2,c/2)\subset \Fol$,
belongs to the set $K(n,c)$. Theorem (\ref{main}) states that
the closure of $K(n,c,\rtt)$ is precisely the set $R_n(c/2,c/2),$
therefore, it is an irreducible component of the space of foliations $\Fol$.

As a consequence of Theorem (\ref{main}), the results in \cite{GML}
and in \cite[Theorem 3.5]{C}, we obtain as a Corollary, a
complete classification of the family of foliations in $K(n,c)$.

\begin{Cor}
  Let $\o\in K(n,c)$ be a foliation, then $K(\o)$ is a complete
  intersection and moreover
  \begin{enumerate}
  \item $\o$ has a meromorphic first integral.
  \item The sets
  $K(n,c,\eta_{pq})=R_n\left(\frac{cp}{p+q},\frac{cq}{p+q}\right),$ so
  that they are irreducible components of the set $\Fol$.
  \item The generic element of $K(n,c,\eta_{pq})$ is structurally stable.
  \end{enumerate}
\end{Cor}

The proof of the Theorem (\ref{main}) uses the ideas introduced in
\cite{CS} and \cite{CMP}, theory of foliations by curves in the
projective plane $\P^2$ and intersection theory, and it will be a
consequence of several lemmas.

\medskip
\subsection{Properties of the Kupka Set}
\label{sec:PK}

From \cite{M, GML} and the references therein, it is known that the
Kupka set is a smooth submanifold of codimension two and has the
\emph{local product structure\/}: For any connected component
$K_0\subset K(\o)$ of the Kupka set, there are an open covering
 ${\mathcal{U}}=\{U_i\}$ of $K_0$, a collection of submersions
$\psi_i:U_i\to \C^2$ such that $K_0\cap U_i=\psi^{-1}_i(0)$, as well a
holomorphic 1--form $\eta=a(x,y)dx-b(x,y)dy$ with a unique singularity
at $(0,0)$, such that the foliation $\o|_{U_i}$ is represented by
$\psi_i^{\ast}\eta$. The one form $\eta$ is called the
\emph{transversal type at the component\/} $K_0$, and it is unique
up to both biholomorphism and multiplication by a never vanishing
holomorphic function.

It is also proved in \cite{GML} that, when the connected component
$K_0$ is compact, and its normal bundle $N_{K_0}$ has non zero
first Chern class, the transversal type may be represented by a
1--form with Taylor expansion $\eta=\eta_{pq}+\cdots$, where
$\eta_{pq}=\tt$ with $p,q$ coprime integer or $p=q=1$.
Furthermore, when the linear part $\eta_{pq}$ belongs to the
Poincar\'e domain, the transversal type is analytically
linearizable. The radial type corresponds to $\eta_{11}=\rtt$. We
denote in this case $R=K(\o)$.

As we will see in the next section, the normal bundle of a compact connected
component of the Kupka set, of a foliation of the projective space
$\P^n,\, n\geq3,$ always has a non zero first Chern class,
and its transversal type has linear part in the Poincar\'e domain.

\medskip

\subsection{Geometric Properties of Kupka Components}
\label{sec:GPK}
The local product property, implies that for any foliation $\o\in
K(n,c)$, it's Kupka set $K=K(\o)$ is subcanonically embedded,
\cite[Theorem 3.4]{CS}. In fact, we have
\begin{equation}
  \label{eq:FN}
    \wedge^2 N_{K}\simeq{\mathcal{O}}_{K}(c),\quad
\mbox{and}\quad\Om^{n-2}_{K}\simeq {\mathcal{O}}_{K}(c-n-1),
\end{equation}
where $\Om^{n-2}_{K}$ denotes the canonical bundle of the Kupka
set $K$.
\medskip

This property has two consequences:

  \begin{description}
  \item[(a)] The first Chern class $c_1(N_{K})\in H^2(K,\Z)$ does
  not vanish. In fact, since $\wedge^2 N_{\K}\simeq
  \mathcal{O}_{K}(c)$ we have that
  $c_1(N_{K})=\j^{\ast}(c\cdot\mathbf{h})$, where
  $\j:K\hookrightarrow\P^n$ denotes the inclusion map, and $\mathbf{h}$
  denotes the class of a hyperplane, i.e. generator of the cohomology ring
  \[ H^{\ast}(\P^n,\Z)\simeq \frac{\Z[\mathbf{h}]}{\mathbf{h}^{n+1}}\]

\item[(b)] By a Serre construction, the normal bundle $N_{K}$ of
  $K$ in $\P^n$, extends to a rank two holomorphic vector bundle
  $V$ of $\P^n$, having a holomorphic section $\si$, vanishing on
  $K(\o)$,  and defining the exact sequence
\begin{equation}\label{exact}
    0\longrightarrow
    {\mathcal{O}}\stackrel{\cdot\si}{\longrightarrow} V
    \longrightarrow {\mathcal{J}}_{K}(c)\rightarrow 0
\end{equation}
where ${\mathcal{J}}_{K}$ denotes the sheaf of ideals of the
Kupka set $K$.
  \end{description}

Since $K$ is smooth, the total Chern class of the bundle $V$ is

\begin{equation}\label{eq:Chern}
    c(V) = 1 + c\cdot \mathbf{h} + deg(K(\o))\cdot \mathbf{h}^2\in
\frac{\Z[\mathbf{h}]}{\mathbf{h}^{n+1}}\simeq H^{\ast}(\P^n,\Z).
\end{equation}
The Kupka set $K$ is a complete intersection if and only if, the
holomorphic vector bundle $V$ splits holomorphically into a direct
sum of holomorphic line bundles.

In what follows, we are going to looking for a condition for the
splitting of the vector bundle $V$, associated to the Kupka set.

Observe that $c_1(V)>0$ and $c_2(V)>0$, and $c_1(N_K)\neq0$, then
the linear part of the transversal type $\eta$ is $\eta_{pq}$ for
some coprime integers $p,q$. Now, the Baum--Bott residues formula
\cite{BB}, implies that
\begin{equation}
  \label{eq:deg}
  deg(K(\o))=\frac{cp}{p+q}\frac{c q}{p+q},
\end{equation}
and the linear part $\eta_{pq}$ belongs to the Poincar\'e domain,
consequently on a neighborhood of each point of the Kupka set, the
foliation looks like the base points of a rational map. Moreover, we
get the decomposition
\[ K(n,c)=\bigcup_{(p,q)\in \mathcal{K}(c)} K(n,c,\eta_{pq}),\]
where $K(n,c,\eta_{pq})$ denotes the set of foliations $\o\in K(n,c)$
with transversal type $\eta_{pq}=px\,dy-qy\,dx$, and we set
\[\mathcal{K}(c)=\left\{ (p,q)\in \N\times\N \,\Big|\,
  \left(\frac{pc}{p+q},\frac{qc}{p+q}\right)\in\N\times\N\right\}\]

The total Chern class of the vector bundle $V$ may be written
\[ c(V)=\left(1+\frac{cp\cdot \mathbf{h}}{p+q}\right)
 \left(1+\frac{qc\cdot\mathbf{h}}{p+q}\right)\]

In particular, if the transversal type is the radial foliation $\rtt$,
we have the following result

\begin{Lemma}\label{degree}
  Let $\o\in K(n,c,\rtt)$ then
  \begin{description}
  \item[(a)] The Chern class of the foliation is even.
  \item[(b)] The total Chern class of the associated vector bundle $V$
  is \[ c(V)=\left(1+\frac{c\cdot\mathbf{h}}{2}\right)^2\]
  \item[(c)] The degree of the Kupka set is \[deg(K)=\frac{c^2}{4}\]
  \end{description}
\end{Lemma}
\noindent\textbf{Proof.} Observe that equations
(\ref{eq:Chern}) and (\ref{eq:deg}) imply
\begin{equation}\label{Chernrad}
   c(V)  =  \left(1 + \frac{c\cdot\mathbf{h}}{2} \right)^2 \quad
  \mbox{and}\quad deg(K(\o))  =  \frac{c^2}{4} = \left(\frac{c}{2} \right)^2,
\end{equation}
therefore, c must be even.$\qed$
\medskip

\subsection{Holomorphic Vector Bundles on the Projective Space}
\label{sec:HVB}
Let us recall some notions of holomorphic vector bundles on the
projective space. Our main reference in this section is \cite{OSS}.

For a rank two holomorphic vector bundle $V$,
which has even first Chern class $c=c_1(V)$, put
\[ V_{norm} := V\left(-\frac{c}{2}\right) \]
and it is normalized in order to have $c_1(V_{norm})=0$.

We say that $V$ is \emph{stable} if $H^0(\P^n,V_{norm})=0$ and $V$
is \emph{semistable} if $H^0(\P^n,V_{norm})\neq0$ and
$H^0(\P^n,V_{norm}(-1))=0.$

Recall that the vector bundle $V$ associated to a Kupka
component never is stable \cite[Corollary 4.3]{CS}. For a foliation
$\o\in K(n,c,\rtt)$ with Kupka set $R$, the total Chern class is
compatible with the semiestability. A Particular case of the next
result was proved and used in \cite{CMP}.

% Recall that the total Chern class of the vector bundle $V$,
% associated to the Kupka set $R$ of a foliation $\o\in K(n,c,\rtt)$
% is compatible with the semistability, moreover it is never stable
% \cite[Corollary 4.3]{CS}. A particular case of the next result was
% proved and used in \cite{CMP}.

\begin{Lemma}\label{KRCI}
Let $\o\in K(n,c,\rtt)$ be a foliation of the projective space
$\P^n$ with Kupka component $R=K(\o)$ and associated vector bundle
$V$. If the holomorphic vector bundle $V$ is semistable and not
stable, then $R$ is a complete intersection.
\end{Lemma}
\noindent{\bf{Proof.}} Let $\tau$ be a non trivial holomorphic
section of the bundle $V_{norm}$, and let $Z$ be the its scheme of
zeros. If $codim(Z)=1$, and $Z_{n-1}$ denotes the union of irreducible
components of dimension $n-1$ (i.e. codimension one), then there is
a homogeneous polynomial $F$ of maximal degree $0<e$ such that
$Z_{n-1}=\{F=0\}$. It follows that $(\tau/F)\in
H^0(\P^n,V_{norm}(-e))$ is a non zero section, contradicting the
semistability of $V$.

If $codim(Z)=2$, we have the exact sequence
\[ 0 \longrightarrow
{\mathcal{O}}\stackrel{\cdot\tau}{\longrightarrow}V_{norm}
\longrightarrow{\mathcal{J}}_Z\longrightarrow0\]

We have that $deg(Z)=c_2(V_{norm})=0$, then $Z$ is empty and
$V_{norm}$ is obtained as an extension of the trivial line bundle by
the trivial line bundle,
hence $V_{norm}\simeq {\mathcal{O}}\oplus {\mathcal{O}}$. Then
\[ V \simeq V_{norm}\left( \frac{c}{2}\right)
\simeq {\mathcal{O}}\left( \frac{c}{2}\right)\oplus
{\mathcal{O}}\left( \frac{c}{2}\right)\qed\]

\medskip
Therefore, the question of when a foliation $\o\in K(n,c,\rtt)$ has a
meromorphic first integral, is reduced to deciding when the rank two
holomorphic vector bundle $V$ associated with the Kupka set $R=K(\o)$
is semistable.

We are going to use the following geometric criteria,
in order to prove the semistability of a rank two vector bundle on
$\P^n$ which belongs to a locally complete intersection.
The proof of the following lemma may be found in \cite[page 186]{OSS}.

\begin{Lemma}\label{CI}
Let $V$ be a rank two holomorphic vector bundle which belongs to the
2 codimensional locally complete intersection $Y\subset \P^n,\quad
n\geq3$ and with $\wedge^2 N_Y = {\mathcal{O}}_Y(c)$. If $c$ is even, then
$V$ is semistable if and only if $c\geq0$ and $Y$ lies in no
hypersurface of degree $d<c/2$.
\end{Lemma}

By \cite[Corollary 4.3]{CS}, the Kupka of a foliation
$\o\in K(n,c,\rtt)$ is contained on an hypersurface of degree $d\leq
c/2$. A combination of the lemmas (\ref{CI}) and (\ref{KRCI})
implies that $\o$ in $K(n,c,\rtt),$ has a meromorphic first integral
if it's Kupka set $R=K(\o)$ no lies in no hypersurface of degree
$d<c/2$. We are going to estimate the degree of an hypersurface $H$
such that $R\subset H$.

% A combination of the lemmas (\ref{CI}) and (\ref{KRCI}) implies that
% any foliation $\o$ in $K(n,c,\rtt),$ has a meromorphic first
% integral if it's Kupka set $R=K(\o)$ no lies in no hypersurface of
% degree $d<c/2$. We are going to estimate the degree of an
% hypersurface $H$ such that $R\subset H$.

\subsection{The Normal Bundle of $R$}
\label{sec:NB}

Given a foliation $\o\in K(n,c,\rtt)$, with Kupka set $R$, we will see
that $R$ can not be contained in the smooth points of an
hypersurface $H$ of degree $d<c/2$.
In order to prove this, we remark the following properties
of the normal bundle $N_R$ of $R$ in $\P^n$.

The first remark is the following: The normal bundle $N_R$ of $R$ in
$\P^n$, is projectively flat, i.e. the projective bundle $\P(N_R)$ is flat.
This property is a consequence of the fact that, the strict transformation
of the foliation obtained by the blowing up along $R$ of a
foliation $\o\in K(n,c,\rtt)$, is transversal to
the exceptional divisor $E=\P(N_R)$. The intersection of this
foliation with $E$, provides the flat structure on the projective
bundle $\P(N_K)$.

On the other hand, for any $\a\in\Z$, the holomorphic vector bundle
$N_R(\a)$ is also projectively flat. Moreover, for any $0<d<c/2$,
the vector bundle $N_R(-d)$ is projectively flat and has positive
first Chern class.

The curvature  matrix
$\widetilde{\Theta}=\Theta(N_R(-d))$ of the twisted normal bundle
$N_R(-d)$ is of the type
\[ \widetilde{\Theta}=\Theta\cdot\mathbf{I}, \]
where $\mathbf{I}\in End(N_R)$ is the identity endomorphism, and
$\Theta$ is a positive 2--form on $R$ \cite[page 7]{K}. Therefore,
the twisted normal bundle $N_R(-d)$ is Griffiths positive \cite{G}, and
hence, it is ample.

We are going to use the following result, which implies the
semistability of the vector bundle $V$, associated to the Kupka
set $R$. The proof may be found in the book \cite[page 133]{FOV}.

\begin{Lemma}\label{TTang} Let $X\subset\P^n$ be a smooth
  subvariety such that the
  twisted normal bundle $N_X(-d)$ is ample. Let $H\subset \P^n$ by a
  hypersurface of degree $d$. Then the set of tangencies
\[ \Gamma=\{ x\in X\cap H_{reg} |T_xX\subset T_xH \}\]
is $0$--dimensional
\end{Lemma}

\medskip
\subsection{Foliations and Curves in Surfaces}
\label{sec:FC}
The main reference through this section, is the book \cite{Br}.

Recall that a holomorphic foliation ${\f}$ by curves on a surface
$M$, may be represented as either
\begin{itemize}
\item By a bundle map $\mathbf{X}:T_{\f}\to TM$.
\item By a holomorphic section $\o$ of the twisted sheaf of one forms
  $\Om^1_M(N_{\f})$.
\end{itemize}
where $T_{\f}$ and $N_{\f}$ are holomorphic line bundles over $M$,
called respectively the \emph{tangent\/} and the \emph{normal\/} bundle of the
foliation. These line bundles are related by the formula
\[ \Om^2_M=T^{\ast}_{\f}\otimes N^{\ast}_{\f}.\]

The above formula, gives for a foliation ${\f}$ on the projective
plane $\P^2$, represented by a section $\o$ of the twisted sheaf
$\Om_{\P^2}^1(c)$, has a normal bundle $N_{\f}={\O}_{\P^2}(c)$ and
its tangent bundle $T_{\f}={\O}_{\P^2}(3-c)$.

Consider a foliation $\o\in{\f}(2,c)$, and assume that there is an
irreducible algebraic curve $C\subset\P^2$ of degree $d$. We are
going find some relations between the degree of the curve $C,$
with the Chern class of the foliation, depending of course, on the
behavior of the curve with respect the foliation. We consider two
cases:
\begin{enumerate}
    \item Non invariant curves.
    \item Invariant curves.
\end{enumerate}

(1) Following \cite[page 23]{Br}, given a foliation $\f$ by curves
on a surface $M$, and a non $\f$--invariant curve $C$, for each
$p\in C$, we define the index $tan({\f},C,p)$, representing the
tangency order of the foliation ${\f}$ with the curve $C$ at $p$ as
follows: let $\{f=0\}$ be a local equation of $C$ around $p$, and
let $\mathbf{X}$ a vector field with isolated singularities and
defining the foliation around $p$. We set
\[ tang({\f},C,p)=dim_{\C}\frac{\O_p}{\langle\mathbf{X}(f),f\rangle}, \]
remark that $\mathbf{X}(f)$ does not vanishes identically on $C$
because it is not ${\f}$--invariant, and so, the ideal
$\langle\mathbf{X}(f),f\rangle$, has finite codimension on $\O_p$ and
$tang({\f},C,p)$ is a finite non negative integer, moreover it is
not zero at points where $C$ is not transverse to ${\f}$, of
course, the singularities of the curve, and the singularities of
the foliations at the curve, are tangency points of the curve $C$ with
the foliation $\f$.

For instance, let $\mathbf{X}=x\partial/\partial x+y\partial/\partial
y$ and $C=\{f=0\}$ a curve through $0\in \C^2.$ Assume that the curve
is singular at $0\in\C^2$, then it
may be written as $f=f_\nu+ f_{\nu+1}+\cdots,\quad \nu\geq2$. In this
case, the index
$tang({\f},C,0)\geq3$, since
\[ \langle \mathbf{X}(f),f\rangle=\langle \nu
f_\nu+\cdots,f_\nu+\cdots\rangle = \langle g,f\rangle, \]
and $g=\nu\cdot f_{\nu+1}+\cdots$.

Now, we define

\[ Tang({\f},C)=\sum_{p\in C}tang({\f},C,p) \]

(2) If the curve $C$ is invariant by the foliation $\o$; we
introduce the following index $Z({\f},C,p)$ for any $p\in C$ as
follows (see \cite[page 24]{Br}). Let $f=0$ be any local equation of
the curve $C$ at $p$, and let $\theta$ be a generating 1--form of
the foliation $\f$ at $p$ . Because $C$ is invariant, we have
\[ g\cdot \theta= h df +f\eta,\quad g,h\in{\O}_p\quad \eta\in\Om^1_p,\]
moreover $(f,h)=1$ in $\O_p$, therefore $h,g$ does not vanish
identically on any branch of $C$. We can define
\[ Z({\f},C,p)=\mbox{ vanishing order of }
\left(\frac{f}{g}\Big|_C\right)\mbox{ at }p .\]
and assuming that $C$ is compact, we define
\[ Z({\f},C)=\sum_{p\in C}Z({\f},C,p). \]

This idex could be non positive for dicritical singularities, for example, let
$\o=\rtt$ be a 1--form generating a foliation on $\C^2$, let
$C=\{f=0\}$ be an invariant curve by the foliation around $p=0$. If $C$ is
smooth (we can take $f=x$), then $Z({\f},C,0)=1$, but it is non
positive in if the curve $C$ is singular.

The result that we need is the following:

\begin{Prop}\label{curves}
  Let $\o\in {\f}(2,c)$ be a foliation on the projective plane
  $\P^2$. Let $C$ be an irreducible algebraic curve of degree $d$ in
  $\P^2$. Then
  \begin{itemize}
  \item If the curve $C$ is invariant by the foliation $\o$ we have
    \begin{equation}
      \label{eq:inv}
     Z({\f},C)=c\cdot d - d^2
    \end{equation}

      \item If the curve $C$ is not invariant by the foliation $\o$.
        \begin{equation}
          \label{eq:notinv}
        Tang({\f},C) = d^2 + (c-3)\cdot d
        \end{equation}
  \end{itemize}
\end{Prop}

\medskip

We can now complete the proof of the Theorem.

\noindent\textbf{Proof of theorem (\ref{main})}.

Let $\o\in K(n,c,\rtt)$, by lemma (\ref{KRCI}), it is sufficient
to prove that the associated vector bundle $V$ is semistable.

If the vector bundle $V$ is not semistable, since $R$ is irreducible,
then there exists an irreducible hypersurface $H$ of degree $d<c/2$
such that $R\subset H$. By the lemma (\ref{TTang}), $R$ is not
contained in the smooth points of the hypersurface $H$.

Take a two plane $\j:\P^2\hookrightarrow \P^n$ in generic position
with respect the foliation and the hypersurface $H$. Consider the
curve $C=\j^{-1}(H)$, and the foliation $\j^{\ast}\o\in
{\f}(2,c).$ This foliation has $\{p_j\}_{j=1,\dots
c^2/4}=\j^{-1}(R)$ radial singularities ($\j\tcap R$), and the
other singularities are Morse type (i.e. there is a local
coordinate system $(x,y)$ such that $\j^{\ast}\o=d(x y)$). Moreover, the
radial singularities are contained in the curve $C$.

If the hypersurface $H$ is not invariant, the same holds for the curve
$C$, moreover, because $C$ is singular at the points $\{p_j\}$ we have
seen that
\[ tang(\j^{\ast}\o,C,p_j)\geq3\quad\mbox{for any}\quad p_j\in \j^{-1}(R).\]
and then, by equation (\ref{eq:notinv}) of proposition (\ref{curves})
we have that $d>c/2$.

If the hypersurface $H$ is invariant, then $Z(\j^{\ast}\o,C,p_j)$ is
non positive, and then, by formula (\ref{eq:inv}) of proposition
(\ref{curves}), its degree $d\geq c>c/2$, therefore we have a
contradiction, therefore, the vector bundle $V$ is
semistable.

By Lemma (\ref{CI}), the holomorphic
vector bundle $V$ associated to $R$ is semistable and by Lemma
(\ref{KRCI}) it is isomorphic to
\[ V\simeq {\mathcal{O}}(c/2)\oplus
{\mathcal{O}}(c/2),\]
Therefore, after the isomorphism between $V$ and
${\mathcal{O}}(c/2)\oplus{\mathcal{O}}(c/2)$, the section $\si=(f,g),$
for some $f,g\in
H^0(\P^n,{\mathcal{O}}_{\P^n}(c/2)).$
Hence, the Kupka set $R$ is a complete intersection and by
Cerveau--Lins theorem \cite{CL}, the foliation has a rational first
integral, and it is defined by the 1--form
\[ \o=f\,dg-g\,df.\]
This finishes the proof of the Theorem \ref{main}$\qed$

\medskip
\begin{Observ}
  Let $M$ be an algebraic manifold with
  $\mbox{dim}_{\C}M\geq3$. Consider the set $K(M,L,\rtt)$
  of those foliations $M$ with normal bundle $L$, (i.e. represented by
  a section $\o$ of the sheaf $\Om^1_M(L)$), and having a compact
  and connected radial Kupka component.

  If the line bundle $L$ is positive, there exists the pair
  a rank two holomorphic vector bundle $V\to M$ with a section $\si$
  inducing the exact sequence
\begin{equation}\label{exact1}
    0\longrightarrow
    {\mathcal{O}}\stackrel{\cdot\si}{\longrightarrow} V
    \longrightarrow {\mathcal{J}}_{K}(L)\rightarrow 0.
\end{equation}

   The vector bundle $V$ has total Chern class
   \begin{equation*}
     c(V)=\left(1+\frac{c_1(L)}{2}\right)^2\in H^{\ast}(M,\Z)
   \end{equation*}
  but in general, the vector bundle $V$, is not semistable, as was proved
  for the projective space.

  In the general case, we have the following conjecture.
\end{Observ}
\begin{Conj}
  Let $\o\in K(M,L,\rtt)$ where $L$ is a very ample line bundle, then
  the associated vector bundle $V$ is projectively flat.
\end{Conj}

This means that any foliation $\o\in
K(M,L,\rtt)$, has a projectively transversal structure. Moreover, if
we consider the section $[\si],$ as a meromorphic section of the
projective bundle $\P(V)$, it has the base locus at $R$, and it defines the
developing map of the transversal structure of the foliation.
\vskip 15pt

\bibliographystyle{amsalpha}

\end{document}